\documentclass[12pt,reqno]{amsart}
\usepackage{enumerate, latexsym, amsmath, amsfonts, amssymb, amsthm, color}
\usepackage{psfrag,graphicx, colordvi}
\newtheorem{theo}{Theorem}

\newtheorem{lem}[theo]{Lemma}

\def\pmod #1{\ ({\rm{mod}}\ #1)}

\def\N{\mathbb N}

\def\Q{\mathbb Q}

\def\l{\left}
\def\r{\right}
\def\bg{\bigg}
\def\({\bg(}
\def\){\bg)}
\def\t{\text}
\def\f{\frac}

\def\ls{\leqslant}

\def\bi{\binom}

\def\eq{\equiv}

\def\Proof{\noindent{\it Proof}}

\theoremstyle{plain}

\theoremstyle{definition}

\theoremstyle{remark}
\newtheorem{remark}{Remark}

\makeatletter \@addtoreset{equation}{section}
\@addtoreset{theo}{section} \makeatother

\def\qed{\hfill \rule{4pt}{7pt}}

\begin{document}
\baselineskip=17pt
\medskip

\hbox{Preprint, {\tt arXiv:1808.04717}}
\medskip

\title
[$q$-Analogues of some series for powers of $\pi$]
{$q$-Analogues of some series for powers of $\pi$}

\author
[Qing-Hu Hou and Zhi-Wei Sun]
{Qing-Hu Hou and Zhi-Wei Sun}

 \address {(Qing-Hu Hou) School of Mathematics, Tianjin University,
  Tianjin 300350, People's Republic of China}

\email{{\tt qh\_hou@tju.edu.cn}
\newline\indent
{\it Homepage}: {\tt http://cam.tju.edu.cn/\lower0.5ex\hbox{\~{}}hou}}

\address {(Zhi-Wei Sun) Department of Mathematics, Nanjing
University, Nanjing 210093, People's Republic of China}

\email{{\tt zwsun@nju.edu.cn}
\newline\indent
{\it Homepage}: {\tt http://maths.nju.edu.cn/\lower0.5ex\hbox{\~{}}zwsun}}


\thanks{The first and second author are supported by the
  National Natural Science Foundation of China (grants 11771330 and 11971222,
  respectively)}

\begin{abstract} We obtain $q$-analogues of several series for powers of $\pi$. For example,
the identity
$$\sum_{k=0}^\infty\f{(-1)^k}{(2k+1)^3}=\f{\pi^3}{32}$$
has the following $q$-analogue:
\begin{equation*}
\sum_{k=0}^\infty(-1)^k\frac{q^{2k}(1+q^{2k+1})}{(1-q^{2k+1})^3}=\frac{(q^2;q^4)_{\infty}^2(q^4;q^4)_{\infty}^6}
{(q;q^2)_{\infty}^4},
\end{equation*}
where $q$ is any complex number with $|q|<1$. We also give $q$-analogues of four new series
for powers of $\pi$ found by the second author.
\end{abstract}

\maketitle

\section{Introduction}

The Riemann zeta function is given by
$$\zeta(s):=\sum_{n=1}^\infty\f1{n^s}\ \quad\t{for}\ \t{Re}(s)>1.$$
Obviously,
$$\sum_{k=0}^\infty\f1{(2k+1)^{s}}=\l(1-\f1{2^s}\r)\zeta(s)\ \ \text{for}\ \t{Re}(s)>1.$$
As Euler proved (cf. \cite[pp.\,231--232]{IR}), for each $m=1,2,3,\ldots$ we have
$$\zeta(2m)=(-1)^{m-1}\f{2^{2m-1}\pi^{2m}}{(2m)!}B_{2m},$$
where the Bernoulli numbers $B_0,B_1,\ldots$ are given by $B_0=1$ and
$$\sum_{k=0}^n\bi{n+1}kB_k=0\ \ \ (n=1,2,3,\ldots).$$

As usual, for $n\in\mathbb N=\{0,1,2,\ldots\}$, the $q$-analogue of $n$ is given by
$$[n]_q = \frac{1-q^n}{1-q}=\sum_{0\ls k<n}q^k.$$
For complex numbers $a$ and $q$ with $|q|<1$, we adopt the standard notation
$$(a;q)_\infty = \prod_{k=0}^{\infty} (1-a q^k).$$
Recently, Z.-W. Sun \cite{S18} obtained the following $q$-analogues of Euler's formulae
$\zeta(2)=\pi^2/6$ and $\zeta(4)=\pi^4/90$:
$$\sum_{k=0}^\infty\f{q^k(1+q^{2k+1})}{(1-q^{2k+1})^2}=\f{(q^2;q^2)_{\infty}^4}
{(q;q^2)_{\infty}^4}$$
and
$$\sum_{k=0}^\infty\f{q^{2k}(1+4q^{2k+1}+q^{4k+2})}{(1-q^{2k+1})^4}
=\f{(q^2;q^2)_{\infty}^8}{(q;q^2)_{\infty}^8},$$
where $q$ is any complex number with $|q|<1$. Note that $\lim_{q\to1}\f{1-q}{1-q^{2k+1}}=\f1{2k+1}$ and also
\begin{equation}\label{1.1}\lim_{q\to1\atop |q|<1}(1-q)\f{(q^2;q^2)_{\infty}^2}{(q;q^2)_{\infty}^2}
=\lim_{q\to1\atop |q|<1}\prod_{n=1}^\infty\f{[2n]_q^2}{[2n-1]_q[2n+1]_q}=\f{\pi}2
\end{equation}
with the help of Wallis' formula
$$\prod_{n=1}^\infty\f{4n^2}{4n^2-1}=\f{\pi}2.$$

Motivated by Sun's work \cite{S18}, A. Goswami \cite{G}
got $q$-analogues of Euler's general formula for $\zeta(2m)$ with $m$ a positive integer,
and M. L. Dawsey and K. Ono \cite{DO} gave further applications.

Let $\chi$ be a Dirichlet character modulo a positive integer $m$. The Dirichlet $L$-function associated with the character $\chi$ is given by
$$L(s,\chi):=\sum_{n=1}^\infty\f{\chi(n)}{n^s}\quad\t{for}\ \Re(s)>1.$$
The Dirichlet beta function is defined by
$$\beta(s)=L\l(s,\l(\f{-4}{\cdot}\r)\r)=\sum_{k=0}^\infty\f{(-1)^k}{(2k+1)^s}\quad\t{for}\ \Re(s)>0,$$
where $(-)$ denotes the Kronecker symbol. As Euler observed,
\begin{equation}\label{1.2}\beta(2n+1)=\f{(-1)^nE_{2n}}{4^{n+1}(2n)!}\pi^{2n+1}
\end{equation}
for all $n=0,1,2,\ldots$ (cf. (3.63) of \cite[p.\,112]{EL}), where $E_0,E_1,E_2,\ldots$ are Euler numbers defined by
$$E_0=1,\ \t{and}\ \sum^n_{k=0\atop 2\mid k}\bi nkE_{n-k}=0\ \ \t{for}\ n=1,2,3,\ldots.$$
In particular,
$$\beta(1)=\f{\pi}4,\ \beta(3)=\f{\pi^3}{32},\ \beta(5)=\f{5\pi^5}{1536}.$$

In view of \eqref{1.1}, we may view Ramanujan's formula
\begin{equation*}
\sum_{k=0}^\infty\f{(-q)^k}{1-q^{2k+1}}=\f{(q^4;q^4)_{\infty}^2}{(q^2;q^4)_{\infty}^2}\ \ (|q|<1)
\end{equation*}
(equivalent to Example (iv) in \cite[p.\,139]{B91}) as a $q$-analogue of
Leibniz's identity $\beta(1)=\pi/4$. Recently, Q.-H. Hou, C. Krattenthaler and Z.-W. Sun \cite{HKS}
obtained the following new $q$-analogue of Leibniz's identity:
\begin{equation*}
\sum_{k=0}^\infty\f{(-1)^kq^{k(k+3)/2}}{1-q^{2k+1}}
=\f{(q^2;q^2)_{\infty}(q^8;q^8)_{\infty}}{(q;q^2)_{\infty}(q^4;q^8)_{\infty}}\quad \ \t{for}\ |q|<1.
\end{equation*}

Motivated by the above work, we seek for a $q$-analogue of the identity
\begin{equation}\label{1.3}\beta(3)=\sum_{k=0}^\infty\f{(-1)^k}{(2k+1)^3}=\f{\pi^3}{32}.
\end{equation}
This leads to our following result.

\begin{theo}\label{Th1.1} For $|q|<1$ we have
\begin{equation}\label{1.4}
\sum_{k=0}^\infty(-1)^k\frac{q^{2k}(1+q^{2k+1})}{(1-q^{2k+1})^3}=\frac{(q^2;q^4)_{\infty}^2(q^4;q^4)_{\infty}^6}
{(q;q^2)_{\infty}^4}.
\end{equation}
\end{theo}

\begin{remark}\label{Rem1.1} \eqref{1.4} is a $q$-analogue of \eqref{1.3} because
\begin{align*}&\lim_{q\to1\atop |q|<1}(1-q)^3\f{(q^2;q^4)_{\infty}^2(q^4;q^4)_{\infty}^6}{(q;q^2)^4}
\\=&\lim_{q\to1\atop |q|<1}(1-q)^2\f{(q^2;q^2)_{\infty}^4}{(q;q^2)_{\infty}^4}
\times\lim_{q\to1\atop |q|<1}\f{1-q^2}{1+q}\cdot\f{(q^4;q^4)_{\infty}^2}{(q^2;q^4)_{\infty}^2}
\\=&\f{\pi^2}4\times\f{\pi}4=\f{\pi^3}{16}
\end{align*}
in view of \eqref{1.1}. \end{remark}

How to give $q$-analogues of \eqref{1.2} for $n=2,3,4,\ldots$? This problem looks sophisticated.

We will show Theorem \ref{Th1.1} in the next section and present more similar results in Section 3.

Recently, Z.-W. Sun \cite{S20} established the following new identities:
\begin{align}\label{Sun1.1}\sum_{k=0}^\infty\f{k(4k-1)\bi{2k}k^3}{(2k-1)^2(-64)^k}=&-\f1{\pi},
\\\label{Sun1.2}\sum_{k=0}^\infty\f{(4k-1)\bi{2k}k^3}{(2k-1)^3(-64)^k}=&\f2{\pi},
\\\label{Sun1.3}\sum_{k=0}^\infty\f{(12k^2-1)\bi{2k}k^3}{(2k-1)^2 256^k}=&-\f2{\pi},
\\\label{Sun1.77} \sum_{k=1}^\infty\f{(3k+1)16^k}{(2k+1)^2k^3\bi{2k}k^3}=&\f{\pi^2-8}2.
\end{align}
They are (1.1), (1.2), (1.3) and (1.77) of Sun \cite{S20} respectively.
In our second theorem we give $q$-analogues of these four identities.

\begin{theo}\label{Th1.2} For $|q|<1$ we have
\begin{gather}\label{S1.1q}
  \sum_{k=0}^\infty (-1)^k q^{k^2} \frac{[2k]_q ([4k]_q-1)}{([2k]_q-1)^2} \cdot\frac{(q;q^2)_k^3}{(q^2;q^2)_k^3} = - \frac{(q;q^2)_\infty (q^3;q^2)_\infty}{(q^2;q^2)_\infty^2},
\\\label{S1.2q}
  \sum_{k=0}^\infty  (-1)^k q^{k^2+2k} \frac{[4k]_q-1}{([2k]_q-1)_q^3} \cdot\frac{(q;q^2)_k^3}{(q^2;q^2)_k^3}=  \frac{(q;q^2)_\infty (q^3;q^2)_\infty}{(q^2;q^2)_\infty^2},
\\\label{S1.3q}
 \sum_{k=0}^\infty \frac{P_k(q)q^{k^2}}{(1-q)^3 ([2k]_q-1)^2}\cdot \frac{(q;q^2)_k^2 (q^2;q^4)_k} {(q^4;q^4)_k^3} = 2q(1+q)
\frac{(q^2;q^4)_\infty (q^6;q^4)_\infty}{ (q^4;q^4)_\infty^2},
\end{gather}
where $P_k(q)$ denotes
\[q^{12k+1}-3q^{10k+2}+3(2q^2-1)q^{8k+1}-(3q^4-1)q^{6k}+3q^{4k+1}-3q^{2k+2}+2q^3-q,
\]
and also
\begin{gather}
\label{S1.77q}
q \sum_{k=0}^\infty \frac{[3k+4]_q}{[2k+3]_q^2} \cdot\frac{(q;q)_k^3 (-q;q)_k} {(q^3;q^2)_k^3} q^{k(k+5)/2} = (1-q)^2
\frac{(q^2;q^2)_\infty^4}{(q;q^2)_\infty^4}-1-q.
\end{gather}

 \end{theo}
We will prove Theorem~1.2 in Section 4.

\section{Proof of Theorem~\ref{Th1.1}}
\label{sec:2}

\begin{lem}\label{Lem2.1} Let $\chi$ be any Dirichlet character. For $|q|<1$ we have
\begin{equation}\label{2.1}\sum_{n=1}^\infty\chi(n)\f{q^{n-1}(1+q^n)}{(1-q^n)^3}=\sum_{m=1}^\infty\(\sum_{d\mid m}\chi\l(\f md\r)d^2\)q^{m-1}.
\end{equation}
\end{lem}
\Proof. For any positive integer $n$, we have
$$\f1{(1-q^n)^3}=\sum_{k=0}^\infty\bi{-3}k(-q^{n})^k=\sum_{k=0}^\infty\bi{k+2}2q^{kn}.$$
Thus
\begin{align*}&\sum_{n=1}^\infty\chi(n)\f{q^{n-1}(1+q^{n})}{(1-q^{n})^3}
\\=&\sum_{n=1}^\infty\chi(n)q^{n-1}\(\sum_{k=0}^\infty\bi{k+2}2q^{kn}+\sum_{k=0}^\infty\bi{k+2}2q^{(k+1)n}\)
\\=&\sum_{n=1}^\infty\chi(n)q^{n-1}\(1+\sum_{k=1}^{\infty}\l(\bi{k+2}2+\bi{k+1}2\r)q^{kn}\)
\\=&\sum_{n=1}^\infty\chi(n)q^{n-1}\sum_{k=0}^{\infty}(k+1)^2q^{kn}
=\sum_{n=1}^\infty\chi(n)\sum_{d=1}^\infty d^2q^{dn-1}
\\=&\sum_{m=1}^\infty \(\sum_{d\mid m}\chi\l(\f md\r)d^2\)q^{m-1}.
\end{align*}
This concludes the proof of \eqref{2.1}. \qed

\medskip
\noindent{\it Proof of Theorem} 1.1. For $n=1,2,3,\ldots$ let $\chi(n)$ be the Kronecker symbol
$(\f{-4}n)$. With the help of Lemma \ref{Lem2.1},
\begin{align*}\sum_{k=0}^\infty(-1)^k\f{q^{2k}(1+q^{2k+1})}{(1-q^{2k+1})^3}
=&\sum_{n=1}^\infty\chi(n)\f{q^{n-1}(1+q^n)}{(1-q^n)^3}
=\sum_{m=1}^\infty \(\sum_{d\mid m}\chi\l(\f md\r)d^2\)q^{m-1}.
\end{align*}
On the other hand, we have
$$\sum_{m=1}^\infty \(\sum_{d\mid m}\l(\f{-4}{m/d}\r)d^2\)q^{m-1}
=\prod_{n=1}^\infty\f{(1-q^{2n})^6(1-q^{4n})^4}{(1-q^n)^4}$$
by Carlitz \cite[4.1]{C} (or \cite[Theorem 2.5]{AAW}). Therefore
\begin{align*}\sum_{k=0}^\infty(-1)^k\f{q^{2k}(1+q^{2k+1})}{(1-q^{2k+1})^3}
=&\prod_{n=1}^\infty\f{(1-q^{2n})^2(1-q^{4n})^4}{(1-q^{2n-1})^4}
\\=&\prod_{n=1}^\infty\f{(1-q^{4n-2})^2(1-q^{4n})^6}{(1-q^{2n-1})^4}
\\=&\f{(q^2;q^4)_{\infty}^2(q^4;q^4)_{\infty}^6}{(q;q^2)_{\infty}^4}.
\end{align*}
This concludes the proof of \eqref{1.4}. \qed

\section{Other results similar to Theorem 1.1}
\label{sec:3}

For any positive integers $d$ and $m$ with $(-1)^md\eq0,1\pmod4$, it is known (cf. \cite{W}) that
$$\sum_{n=1}^\infty\l(\f{(-1)^md}n\r)\f1{n^m}\in\f{\pi^m}{\sqrt d}\Q,$$
where $\Q$ is the field of rational numbers.

 Let $q$ be any complex number with $|q|<1$.
By Carlitz \cite[(2.1)]{C} (or \cite[Theorem 2.3]{AAW}),
$$\sum_{m=1}^\infty\(\sum_{d\mid m}\l(\f{-3}{m/d}\r)d^2\)q^{m-1}=\prod_{n=1}^\infty\f{(1-q^{3n})^9}{(1-q^n)^3}.$$
Combining this with Lemma \ref{Lem2.1} we obtain that
\begin{equation}\label{3.1}
\sum_{n=0}^\infty\l(\f n3\r)\f{q^{n-1}(1+q^n)}{(1-q^n)^3}=\f{(q^3;q^3)_{\infty}^9}{(q;q)_{\infty}^3},
\end{equation}
which is a $q$-analogue of the identity
\begin{equation}\label{3.2}\sum_{n=1}^\infty\l(\f n3\r)\f1{n^3}=\f{4\pi^3}{81\sqrt3}.
\end{equation}

Similar to Lemma \ref{Lem2.1}, for any Dirichlet character $\chi$, we have
\begin{equation}\label{3.3}\sum_{n=1}^\infty\f{\chi(n)q^{n-1}}{(1-q^n)^2}
=\sum_{n=1}^\infty\chi(n)\sum_{k=0}^\infty(k+1)q^{(k+1)n-1}=\sum_{m=1}^\infty\(\sum_{d\mid m}\chi\l(\f md\r)d\)q^{m-1}.
\end{equation}
Combining this with Ramanujan's identity
$$\sum_{m=1}^\infty\(\sum_{d\mid m}\l(\f 5{m/d}\r)d\)q^{m-1}=\prod_{n=1}^\infty\f{(1-q^{5n})^5}{1-q^n},$$
we recover Ramanujan's result
\begin{equation}\label{3.4}\sum_{n=1}^\infty\l(\f n5\r)\f{q^{n-1}}{(1-q^n)^2}=\f{(q^5;q^5)_{\infty}^5}{(q;q)_{\infty}}
\end{equation}
(cf. \cite[p.\,107]{B06}),
which can be viewed as a $q$-analogue of the identity
\begin{equation}\label{3.5}\sum_{n=1}^\infty\l(\f n5\r)\f1{n^2}=\f{4\pi^2}{25\sqrt5}.
\end{equation}
Ramanujan used \eqref{3.4} to deduce his famous congruence $p(5n+4)\eq0\pmod 5$, where $n$ is any nonnegative integer and $p(\cdot)$ is the well-known partition function.

By \cite[Theorems 3.5 and 3.7]{AAW},
$$\sum_{m=1}^\infty\(\sum_{d\mid m}\l(\f{8}{m/d}\r)d\)q^{m-1}=\prod_{n=1}^\infty\f{(1-q^{2n})^3(1-q^{4n})(1-q^{8n})^2}{(1-q^n)^2}$$
and
$$\sum_{m=1}^\infty\(\sum_{d\mid m}\l(\f{12}{m/d}\r)d\)q^{m-1}=\prod_{n=1}^\infty\f{(1-q^{2n})^2(1-q^{3n})^2(1-q^{4n})(1-q^{12n})}{(1-q^n)^2}.$$
Combining this with \eqref{3.3} we obtain
\begin{equation}\label{3.6}\begin{aligned}
\sum_{k=0}^\infty(-1)^{k(k+1)/2}\f{q^{2k}}{(1-q^{2k+1})^2}
=\f{(q^2;q^2)_{\infty}^3(q^4;q^4)_{\infty}(q^8;q^8)_{\infty}^2}{(q;q)_{\infty}^2}
\end{aligned}\end{equation}
and
\begin{equation}\label{3.7}\begin{aligned}
&\sum_{k=0}^\infty\l(\f{q^{6k}}{(1-q^{6k+1})^2}-\f{q^{6k+4}}{(1-q^{6k+5})^2}\r)
\\=&\f{(q^2;q^2)_{\infty}^2(q^3;q^3)_{\infty}^2(q^4;q^4)_{\infty}(q^{12};q^{12})_{\infty}}{(q;q)_{\infty}^2},
\end{aligned}
\end{equation}
which are $q$-analogues of the identities
\begin{equation}\label{3.8}\sum_{k=0}^\infty\f{(-1)^{k(k+1)/2}}{(2k+1)^2}=\f{\pi^2}{8\sqrt2}
\ \t{and}\ \sum_{k=0}^\infty\l(\f 3{2k+1}\r)\f1{(2k+1)^2}=\f{\pi^2}{6\sqrt3}.
\end{equation}

We are unable to find $q$-analogues of many identities including the following ones£º
\begin{align*}\sum_{k=0}^\infty\f{(-1)^{k(k+3)/2}}{(2k+1)^3}&=\f{3\pi^3}{64\sqrt2},
\\
\sum_{k=0}^\infty\f{(-1)^{k(k+1)/2}}{(2k+1)^4}&=\f{11\pi^4}{768\sqrt2},
\\ \sum_{k=0}^\infty\l(\f3{2k+1}\r)\f1{(2k+1)^4}&=\f{23\pi^4}{1296\sqrt3},
\\\sum_{n=1}^\infty\l(\f n3\r)\f1{n^5}&=\f{4\pi^5}{729\sqrt3}.
\end{align*}

\section{Proof of Theorem 1.2}
\label{sec:4}

\medskip
\noindent{\it Proof of Theorem 1.2}.
For $k\in\N$ let
\[
  a_k(q) = (-1)^k q^{k^2} [4k+1]_q \frac{(q;q^2)_k^3}{(q^2;q^2)_k^3}.
\]
It is known that (see, e.g., \cite[Eq.~(1.5)]{Guo20})
\[
\sum_{k=0}^\infty a_k(q) = \frac{(q;q^2)_\infty (q^3;q^2)_\infty}{(q^2;q^2)_\infty^2},
\]
which is the $q$-analogue of Bauer's formula
\[
\sum_{k=0}^\infty (-1)^k (4k+1) \frac{{2k \choose k}^3}{64^k} = \frac{2}{\pi}.
\]

Let
  \[
  b_k(q) = (-1)^k q^{k^2} \frac{[2k]_q ([4k]_q-1)}{([2k]_q-1)^2}\cdot \frac{(q;q^2)_k^3}{(q^2;q^2)_k^3},
  \]
be the summand on the left hand side of \eqref{S1.1q}.
It is easy to verify that
  \[
  a_k(q) + b_k(q) = \Delta_k\left( (-1)^{k+1} q^{k^2-1} \frac{(1-q^{2k})^3}{(1-q^{2k-1})^2(1-q)} \cdot\frac{(q;q^2)_k^3}{(q^2;q^2)_k^3} \right)
  \]
for all $k\in\N$, where $\Delta_k$ is the difference operator defined by
\[
\Delta_k f(k) = f(k+1) - f(k).
\]
Hence
\[
\sum_{k=0}^\infty b_k(q) = - \sum_{k=0}^\infty a_k(q) = - \frac{(q;q^2)_\infty (q^3;q^2)_\infty}{(q^2;q^2)_\infty^2}.
\]
This proves \eqref{S1.1q}.

 Similarly, we set
  \[
  c_k(q) = (-1)^k q^{k^2+2k} \frac{[4k]_q-1}{([2k]_q-1)^3}\cdot \frac{(q;q^2)_k^3}{(q^2;q^2)_k^3}
  \]
  be the summand on the left hand side of \eqref{S1.2q} and note that
  \[
  - a_k(q) +c_k(q) = \Delta_k\left( (-1)^{k} q^{k^2-2} \frac{(1-q^{2k})^3 (q^2-q^{2k})}{(1-q^{2k-1})^3(1-q)}\cdot \frac{(q;q^2)_k^3}{(q^2;q^2)_k^3} \right),
  \]
for all $k\in\N$. Hence
  \[
  \sum_{k=0}^\infty c_k(q) = \sum_{k=0}^\infty a_k(q)= \frac{(q;q^2)_\infty (q^3;q^2)_\infty}{(q^2;q^2)_\infty^2}.
  \]
This proves \eqref{S1.2q}.

For $k\in\N$ let
\[
s_k(q) = \frac{P_k(q)q^{k^2+1}}{(1-q)^3 ([2k]_q-1)^2}\cdot \frac{(q;q^2)_k^2 (q^2;q^4)_k} {(q^4;q^4)_k^3}
\]
and
\[
t_k(q) = [6k+1]_q \frac{(q;q^2)_k^2 (q^2;q^4)_k}{(q^4;q^4)_k^3} q^{k^2}.
\]
Then
\[
s_k(q) -2q^2t_k(q)= \Delta_k \left(
 \frac{(1-q^{4k})^3 q^{k^2}}{(1-q^{2k-1})^2(1-q)}
 \cdot\frac{(q;q^2)_k^2 (q^2;q^4)_k} {(q^4;q^4)_k^3}
\right).
\]
for all $k\in\N$. Therefore,
\[
\sum_{k=0}^\infty s_k(q) = 2q^2 \sum_{k=0}^\infty t_k(q)
= 2q^2(1+q)\frac{(q^2;q^4)_\infty (q^6;q^4)_\infty}{(q^4;q^4)_\infty^2}
\]
with the aid of \cite[Theorem 1.1]{GL}.
This proves \eqref{S1.3q}.

Finally, let us consider \eqref{S1.77q}. It is easy to verify that
\begin{multline*}
  \frac{q^{2k+1}[3k+4]_q}{[2k+3]_q^2}\cdot \frac{(q;q)_k^3 (-q;q)_k} {(q^3;q^2)_k^3} q^{k(k+1)/2}
-  [3k+2]_q \frac{(q;q)_k^3 (-q;q)_k} {(q^3;q^2)_k^3} q^{k(k+1)/2} \\
  = \Delta_k \left( \frac{(1+q^{k+1})(1-q^{2k+1})}{1-q}\cdot \frac{(q;q)_k^3 (-q;q)_k} {(q^3;q^2)_k^3} q^{k(k+1)/2} \right)
\end{multline*}
for all $k\in\N$. Therefore,
$$
  q \sum_{k=0}^\infty \frac{[3k+4]_q}{[2k+3]_q^2}\cdot \frac{(q;q)_k^3 (-q;q)_k} {(q^3;q^2)_k^3} q^{k(k+5)/2}
  - \sum_{k=0}^\infty [3k+2]_q \frac{(q;q)_k^3 (-q;q)_k} {(q^3;q^2)_k^3} q^{k(k+1)/2}$$
  coincides with $-1-q$.
By  \cite[(1.9)]{HKS},
\[
\sum_{k=0}^\infty [3k+2]_q \frac{(q;q)_k^3 (-q;q)_k} {(q^3;q^2)_k^3} q^{k(k+1)/2}
 = (1-q)^2
\frac{(q^2;q^2)_\infty^4}{(q;q^2)_\infty^4}.
\]
So we have \eqref{S1.77q}.

 In view of the above, we have completed the proof of Theorem 1.2. \qed

\end{document}